\documentclass{article} 

\title{\sf The camel-banana problem}
\author{Michiel de Bondt}
\date{June 12, 1996\footnote{The proof of Lemma B has been replaced by a shorter proof on February 10, 2024.}}

\usepackage{latexsym}
\usepackage{epsf}

\begin{document}

\maketitle

\begin{quotation}
\begin{sf} \small
A camel can carry one banana at a time on its back. It is on a diet and 
therefore can only have one banana at a time in its stomach. As soon as it 
has eaten a banana it walks a mile and then it needs a new banana (in 
order to be able to continue its itinerary).  

Let there be a stock of $N$ bananas at the border of the desert. How far 
can the camel penetrate into the desert, starting at this point? (Of 
course it can form new stocks with transported bananas.)\/\end{sf}%
\footnote{This was problem number 10 of the ``Universitaire Wiskunde
Competitie'' of 1992-93, a Dutch mathematical competition for students.
A partial solution then was given by the author. Here he presents a full
solution.}
\end{quotation}

\bigskip \noindent
It is easy to see that the camel can penetrate $1$ mile into the
desert with $1$ banana, $2$ miles with $2$ bananas and $2\frac13$
miles with $3$ bananas.  But how far can the camel get with
$73083734$ bananas? We shall give the answer in the following pages.

We can assume that the camel will use only one line of the desert. It is 
natural to take a real co-ordinate along this line, with the starting 
point at the origin and a mile as the unit length. Positions along the 
line will be given by their co-ordinate. 
An `Ultra Wise Camel' might apply the following strategy, since it is 
optimal, as we shall prove. 

\bigskip \noindent
UWC-STRATEGY: Repeat the following until there are $3$ bananas left: 
\begin{enumerate}

\item[(1)] If the number of bananas left is {\em even}\/: eat a banana, carry 
          another banana as far as possible and return to the banana(s) 
          closest to the origin.

\item[(2)] If the number of bananas left is {\em odd}\/: eat a banana, carry 
          another banana to the furthest banana, then carry both bananas 
          (each in turn) as far as possible, and finally return to the 
          banana(s) closest to the origin.  

\end{enumerate}
When there are $2$ or $3$ bananas left, do the above without returning to 
the banana(s) closest to the origin. 

\begin{figure}[!t]
\begin{center}
\epsfbox{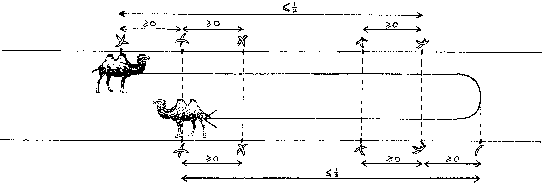}
\caption{a UWC-strategy step for an even number of bananas.}
\end{center}
\end{figure}

\begin{figure}[!t]
\begin{center}
\epsfbox{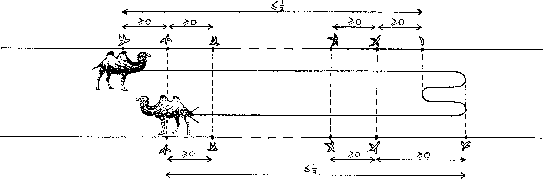}
\caption{a UWC-strategy step for an odd number of bananas.}
\end{center}
\end{figure}

\begin{figure}[!t]
\begin{center}
\epsfbox{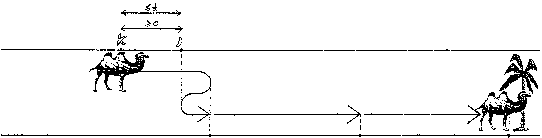}
\caption{the last 3 steps of the UWC-strategy.}
\end{center}
\end{figure}

Figure 3 shows the last $3$ steps of 
this strategy. Steps (1) and (2) are illustrated in Figures 
1 and 2. The figures show that the following properties holding at the
beginning are not affected by applying the UWC-strategy:
\begin{enumerate}

\item[(a)] During meal-time, the camel is with the banana(s) closest to the
          origin and there is an interval of length $\frac{1}{2}$ containing 
          all bananas. 

\item[(b)] Just before step (1), the bananas can be divided into {\em 
          pairs}\/, such that the bananas of each pair are at the same 
          position (but different pairs may also be at the same position). 
          Figure 1 makes clear that the camel can transform the hindmost 
          pair into one furthest banana (with one of the two bananas being 
          eaten and the other being transported) and leave all other 
          pairs where they are in step (1).  

\item[(c)] Just before step (2), the bananas can be divided 
          into pairs in the same way as in (b), except one furthest 
          banana. Figure 2 makes clear that the camel can transform the 
          hindmost pair and the sole furthest banana into a new furthest
          pair and leave all other pairs where they are in step (2).

\end{enumerate}

\bigskip \noindent
Let $c(N)$ denote the distance the camel can penetrate into the desert, when
applying the UWC-strategy.

\bigskip \noindent
EXAMPLES: We illustrate the UWC-strategy for the cases of $N = 3$, $4$, 
$5$, $6$, $7$ bananas by giving the positions of the bananas at each step 
just before the camel starts eating.  

\begin{center}
\begin{tabular}{lp{4in}}
$N=3$: & $(0, 0, 0)$, $(\frac{1}{3}, \frac{1}{3})$, $(1\frac{1}{3})$ 
         $\Longrightarrow c(3) = 2\frac{1}{3}$ \\ 
$N=4$: & $(0, 0, 0, 0)$, $(0, 0, \frac{1}{2})$, $(\frac{2}{3}, \frac{2}{3})$, 
         $(1\frac{2}{3})$ $\Longrightarrow c(4) = 2\frac{2}{3}$ \\ 
$N=5$: & $(0, 0, 0, 0, 0)$, $(0, 0, \frac{1}{4}, \frac{1}{4})$,
         $(\frac{1}{4}, \frac{1}{4}, \frac{5}{8})$, 
         $(\frac{5}{6}, \frac{5}{6})$, $(1\frac{5}{6})$
         $\Longrightarrow c(5) = 2\frac{5}{6}$ \\
$N=6$: & $ (0, 0, 0, 0, 0, 0)$, $(0, 0, 0, 0, \frac{1}{2})$, 
         $(0, 0, \frac{1}{2}, \frac{1}{2})$, 
         $(\frac{1}{2}, \frac{1}{2}, \frac{3}{4})$,
         $(1, 1)$, $(2)$ $\Longrightarrow c(6) = 3$ \\
$N=7$: & $(0, 0, 0, 0, 0, 0, 0)$,
         $(0, 0, 0, 0, \frac{1}{4}, \frac{1}{4})$, 
         $(0, 0, \frac{1}{4}, \frac{1}{4}, \frac{1}{2})$,
         $(\frac{1}{4}, \frac{1}{4}, \frac{9}{16}, \frac{9}{16})$, 
         $(\frac{9}{16}, \frac{9}{16}, \frac{29}{32})$, 
         $(1\frac{1}{8}, 1\frac{1}{8})$, $(2\frac{1}{8})$ 
         $\Longrightarrow c(7) = 3\frac{1}{8}$
\end{tabular}
\end{center}

\bigskip \noindent
Suppose that the camel applies the UWC-strategy. Since pairs of
bananas as described in properties (b) and (c) 
are destructed only during meal-time, each position of a pair is an {\em
eating position}\/: a position where a banana is eaten during the strategy.  

Let $e_i$ denote the eating position where the camel eats its $(N+1-i)$-th 
banana in the UWC-strategy. Then
$$
0 = e_{N} \le e_{N-1} \le \ldots \le e_{2} \le e_{1} = c(N)-1
$$ 
When there are $2n$ bananas left just before a strategy step, the $n$ pairs 
of bananas are lying at $e_{2n}$, $e_{2n-1}$, $e_{2n-2}$, $\ldots$, 
$e_{n+1}$. Inspection of what happens at the next two steps of the 
UWC-strategy reveals that 
\begin{equation}
e_n = \frac{1}{4}e_{2n-2} + \frac{1}{2} e_{2n-1} + \frac{1}{4} e_{2n}
      + \frac{1}{2} \quad (n \ge 2)
\label{nw1}\end{equation}
Note that for $n=2$, this equality simplifies to
$$
e_2 = \frac{2}{3}e_3 + \frac{1}{3}e_4 + \frac{2}{3}
$$ 
This implies that for each $n \ge 2$, the number $e_2$ can be expressed as 
an affinely linear function of $e_{2n}$, $e_{2n-1}$, $e_{2n-2}$, $\ldots$, 
$e_{n+1}$ with coefficients depending only on $n$ and not on the original 
number $N$ of bananas at the starting position. 
The same
holds for $c(N) = e_2 + 2$, in other words, for each $n \ge 2$ with
$2n \le N$ there are constants $\lambda_{n,i}$ ($i = 0$, $1$, $\ldots$, $n$) 
such that
$$ 
c(N) = \lambda_{n,n} e_{2n} + \lambda_{n,n-1} e_{2n-1} + \ldots +
\lambda_{n,1} e_{n+1} + \lambda_{n,0}
$$
Note also that $\sum_{i=1}^{n} \lambda_{n,i} = 1$ for all $n$.
Indeed, this relation holds for $n = 2$, and it is preserved during all
subsequent substitutions. We shall use these observations to prove:

\paragraph{Theorem I.}
If $n \ge 1$, then $c$ satisfies the following equalities
\begin{eqnarray*}
c(2n+1) &=& \frac{1}{2}c(2n) + \frac{1}{2}c(2n+2) \\
c(2n) &=& \frac{1}{2}c(n) + \frac{1}{2}c(n+1) + \frac{1}{2}
\end{eqnarray*}

\paragraph{Proof:}
The examples above show that the first equality holds for $n = 1$,
so let $n > 1$.
Suppose that there are $2n$ bananas left and the camel is about to eat its 
next banana, while applying the UWC-strategy for $N$ bananas 
($N \ge 2n$). Then: 
\begin{itemize}

\item if $N = 2n$ then there are $2n$ bananas at $0$, so 
      $e_{2n}=\ldots=e_{n+1}=0$, and $c(2n) = \lambda_{n,0}$,

\item if $N = 2n+1$ then there are $2n-2$ bananas at $0$ and $2$ bananas 
      at $\frac{1}{4}$, so $e_{2n}=\ldots=e_{n+2}=0$,  
      $e_{n+1}=\frac{1}{4}$, and $c(2n+1) =
      \frac{1}{4}\lambda_{n,1} + \lambda_{n,0}$, 

\item if $N = 2n+2$ then there are $2n-2$ bananas at $0$ and $2$ bananas
      at $\frac{1}{2}$, so $e_{2n}=\ldots=e_{n+2}=0$, 
      $e_{n+1}=\frac{1}{2}$, and $c(2n+2) =
      \frac{1}{2}\lambda_{n,1} + \lambda_{n,0}$. 

\end{itemize}
The first equality of the theorem thus is proved. 
The proof of the second equality is slightly more complicated. We 
split it into two cases, $n = 2m$ and $n = 2m+1$: 
\begin{eqnarray*}
c(4m) &=& \frac{1}{2}c(2m) + \frac{1}{2}c(2m+1) + \frac{1}{2} \\
c(4m+2) &=& \frac{1}{2}c(2m+1) + \frac{1}{2}c(2m+2) + \frac{1}{2}
\end{eqnarray*}
The above examples show that these equalities hold for $m = 1$, 
so let $m > 1$. 
We apply the UWC-strategy for $N=4m$ and $N=4m+2$ separately, and
consider in each case the
situation when there are $2m$ bananas left.
\begin{itemize}

\item If $N = 4m$ then there are $2m-2$ bananas at $\frac{1}{2}$ and $2$ 
      bananas at $\frac{5}{8}$, so $e_{2m}=\ldots=e_{m+2}=\frac{1}{2}$ and 
      $e_{m+1}=\frac{5}{8}$. It follows that
      \begin{eqnarray*}
      c(4m) &=& \frac{1}{2} \sum_{i=2}^m \lambda_{m,i} + \frac{5}{8}
      \lambda_{m,1} + \lambda_{m,0} = \frac{1}{2} +
      \frac{1}{8}\lambda_{m,1} + \lambda_{m,0} \\
      &=& \frac{1}{2}c(2m) + \frac{1}{2}c(2m+1) + \frac{1}{2}    
      \end{eqnarray*}

\item If $N = 4m+2$ then there are $2m-2$ bananas at $\frac{1}{2}$ and $2$ 
      bananas at $\frac{7}{8}$, so $e_{2m}=\ldots=e_{m+2}=\frac{1}{2}$ and 
      $e_{m+1}=\frac{7}{8}$.  It follows that
      \begin{eqnarray*}
      c(4m+2) & = & \frac{1}{2} \sum_{i=2}^m \lambda_{m,i} + \frac{7}{8}
      \lambda_{m,1} + \lambda_{m,0} = \frac{1}{2} +
      \frac{3}{8}\lambda_{m,1} + \lambda_{m,0} \\
      & = & \frac{1}{2}c(2m + 1) +  \frac{1}{2}c(2m+2) + \frac{1}{2} 
      \end{eqnarray*}

\end{itemize}
This completes the proof of Theorem I. \hfill $\Box$

\bigskip \noindent
Although the equalities of theorem I and the initial value $c(1) = 1$ 
completely determine the numbers $c(n)$ for all $n$, to compute, say,
$c(73083734)$, it is more convenient to have an explicit formula. We
now shall derive one.  

Define 
\begin{equation}
F(n,m) = \frac{^2\!\log m}{2} + g(\frac{n-1}{m}) 
          + \frac{(n-1) \pmod 2}{6m^2} \quad (n,m>0)
\label{nw2}
\end{equation}
where
$$
g(x) = -\frac{1}{6}x^2 + x + 1
$$
Then it is a trivial exercise to verify the following equalities:
\begin{eqnarray}
F(2m,2m) &=& F(2m,m) \label{eq:1} \\
F(2n+1,m) &=& \frac{1}{2}F(2n,m)+\frac{1}{2}F(2n+2,m) \label{eq:2} \\
F(2n,2m) &=& \frac{1}{2}F(n,m)+\frac{1}{2}F(n+1,m)+\frac{1}{2} \label{eq:3}
\end{eqnarray}
Define
$$
h(n) = 2^{\lfloor{^2\!\log n}\rfloor} \quad (n \ge 1)
$$
so $h(n)$ is the largest power of $2$ such that $h(n) \le n$.

\paragraph{Theorem II.}
$$
c(n) = F(n, h(n)) \quad (n \ge 1)
$$

\paragraph{Proof:}
The case $n=1$ is trivial, so let $n > 1$. It is sufficient to check the 
equalities corresponding to those of theorem I: 
\begin{eqnarray*}
F(2n+1, h(2n+1)) &=& \frac{1}{2}F(2n, h(2n)) 
                      + \frac{1}{2}F(2n+2, h(2n+2)) \quad (n \ge 1) \\
F(2n, h(2n)) &=& \frac{1}{2}F(n, h(n)) + \frac{1}{2}F(n+1, h(n+1)) 
                  + \frac{1}{2} \quad (n \ge 1)
\end{eqnarray*}
To do so, we consider two cases: 
\begin{itemize}

\item $n+1 = 2^k$ for some integer $k \ge 1$. \\
      Then $h(n) = 2^{k-1}$, $h(n+1) = h(2n)= h(2n+1) = 2^k$ 
      and $h(2n+2) = 2^{k+1}$. By
      (\ref{eq:1}), the two equalities then reduce to
      \begin{eqnarray*}
      F(2n+1, n+1) &=& \frac{1}{2}F(2n,n+1) 
                       + \frac{1}{2}F(2n+2,2n+2)  \\
                   &=& \frac{1}{2}F(2n,n+1)  + \frac{1}{2}F(2n+2,n+1)  \\
      F(2n,n+1)    &=& \frac{1}{2}F(n,\frac{n+1}{2}) + \frac{1}{2}F(n+1,n+1) 
                       + \frac{1}{2} \\
                   &=& \frac{1}{2}F(2n,\frac{n+1}{2}) 
                       + \frac{1}{2}F(n+1,\frac{n+1}{2})  \\
      \end{eqnarray*}
      which are true by (\ref{eq:2}) and (\ref{eq:3}). 

\item $2^{k} < n+1 < 2^{k+1}$ for some integer $k \ge 1$. \\
      Then $h(n) = h(n+1) = 2^k$, $h(2n) = h(2n+1) = h(2n+2) = 2^{k+1}$ and 
      the two equalities can be written as
      \begin{eqnarray*}
      F(2n+1, 2^{k+1}) &=& \frac{1}{2}F(2n,2^{k+1}) 
                           + \frac{1}{2}F(2n+2,2^{k+1}) \quad (n \ge 1) \\
      F(2n,2^{k+1})    &=& \frac{1}{2}F(n,2^k) + \frac{1}{2}F(n+1,2^k) 
                           + \frac{1}{2} \quad (n \ge 1)
      \end{eqnarray*}  
      which are also true by (\ref{eq:2}) and (\ref{eq:3}).  \hfill $\Box$

\end{itemize}

\bigskip \noindent
Now $c(73083734) = 14\frac{1003590240076691}{1125899906842624}$ can
be easily verified using theorem II and formula (\ref{nw2}). But we still
have to show that the camel cannot reach further than this with $73083734$
bananas.

Before proving the optimality of the UWC-strategy, we first
derive two more properties of the UWC-strategy. Let
$$
s_{n} = e_{1} + 2\sum_{k=2}^n e_{k} 
$$
(which by definition also includes $s_{1} = e_{1}$). 

\paragraph{Lemma A.}
$$
\begin{array}{ll} 
s_{n} = N-1                        & \quad (n \ge N/2) \\
s_{n} = \frac{1}{4}s_{2n-1} + 
\frac{1}{4}s_{2n} + \frac{2n-1}{2} & \quad (n \le N/2)
\end{array}
$$

\paragraph{Proof:}
Suppose that the camel has just eaten 
its last banana in the UWC-strategy.
It has walked $N-1$ miles in total, of which, say, $f$ miles in forward 
direction, and $b$
miles backward. Then $f + b = N-1$ and $f - b = e_{1}$, so $f =
(N-1+e_1)/2$. 
Since the camel always had a banana on its back when it
walked forward, the sum $\sum_{n=1}^N e_n$ of the positions of the banana
skins equals $(N-1+e_1)/2$. Since $e_n=0$ if $n \ge N/2 + 1$,
it follows that $s_{n} = N-1$ for all $n \ge N/2$. \\
To prove the remaining part of the lemma, we note that by equation
(\ref{nw1}) 
\begin{eqnarray*}
2e_k &=& \frac12 e_{2k-2} + e_{2k-1} + \frac12 e_{2k} + 1\\
     &=& \frac14(2e_{2k-2} + 2e_{2k-1}) + \frac14(2e_{2k-1} + 2e_{2k}) + 1
\end{eqnarray*}
for $k = 2$, $3$, $\ldots$, $n$. Since also 
$$
e_1 = \frac12e_1 + \frac12e_2 + \frac12 =
\frac14 s_1 + \frac14 s_2 + \frac12
$$
the desired result follows by addition. \hfill $\Box$

\bigskip \noindent
Now let $S$ be any optimal strategy for $N$ bananas.
Let $e_1'$, $e_2'$, $\ldots$, $e_N'$ 
be the eating positions of the $N$ bananas
(which are all used), numbered in such a way that
$$
0 = e_{N}' \le e_{N-1}' \le \ldots \le e_{2}' \le e_{1}' 
$$
Note that we do not suppose that the bananas are eaten in the order 
$e_{N}'$, $e_{N-1}'$, $\ldots$, $e_{1}' $, although we
may (and shall) assume that the first banana is eaten at
$e_{N}'$, and the last one at $e_{1}'$. In between, however, the
order might be different. 

As above, we define 
$$
s_{n}' = e_{1}' + 2\sum_{k=2}^n e_{k}' 
$$
and as in the proof of Lemma A, we suppose that the camel has just eaten
its last banana.  Again, it has walked $N-1$ miles, 
of which $(N-1+e_1')/2$ miles in forward direction. Since it
always had at most one banana on its back when it walked forward, the
sum $\sum_{n=1}^N e_n'$ of the positions of the banana skins is at
most $(N-1+e_1')/2$. This implies that $s_{N}' \le N-1$, and since
$s_{n}' \le s_{N}'$ for all $n$, we also have $s_{n}' \le N-1$ for
all $n$. \\ 
As an analogue of Lemma A, we formulate:

\paragraph{Lemma B.}
$$
\begin{array}{ll}
s_{n}' \le N-1 & \quad (n \ge N/2) \\
s_{n}' \le \frac{1}{4}s_{2n-1}' + \frac{1}{4}s_{2n}' + \frac{2n-1}{2} 
               & \quad (n \le N/2)
\end{array}
$$

\bigskip \noindent
The second inequality of Lemma B is somewhat more difficult to prove. 
Before giving a proof of lemma B, we note that lemmas A
and B together yield $s_{n}' \le s_{n}$ for all $n$ (start with $n = N$ and
proceed in descending order), so in particular 
$e_1' + 1 = s_1' + 1 \le s_1 + 1 = e_1 + 1 = c(N)$. Since $S$ is assumed to
be an optimal strategy, the camel reaches as far as $e_1' + 1 \ge c(N)$ in
$S$. So $e_1' + 1 = c(N)$, which amounts to the optimality of the 
UWC-strategy.

\paragraph{Theorem III.}
The UWC-strategy is optimal.

\paragraph{Proof:}
We are done when we shall have proved lemma B.

Let $\beta = e'_{2n} + \frac12$. During the $\frac12$ mile before eating a banana on position $e'_i$, the camel can maximize the portion of its walk that is done in the interval $[\beta,\infty)$ by walking from position $e'_i + \frac12$ to position $e'_i$. It walks at most $\max\{e'_i+\frac12,\beta\} - \max\{e'_i,\beta\}$ miles in $[\beta,\infty)$. During the $\frac12$ mile after eating a banana on position $e'_i$, the camel can maximize the portion of its walk that is done in the interval $[\beta,\infty)$ by walking from position $e'_i$ to position $e'_i + \frac12$. Again, it walks at most $\max\{e'_i+\frac12,\beta\} - \max\{e'_i,\beta\}$ miles in $[\beta,\infty)$.

For $2 \le i \le n$, we estimate $\max\{e'_i+\frac12,\beta\} - \max\{e'_i,\beta\}$ by $\frac12$. For $n+1 \le i \le 2n$, we estimate $\max\{e'_i+\frac12,\beta\} - \max\{e'_i,\beta\}$ by $\max\{e'_i+\frac12,\beta\} - \beta = e'_i + \frac12 - \beta$. So before eating its last banana on position $e'_1$, the camel walks
$$
\alpha \le (e'_1 + \frac12 - \max\{e'_1,\beta\}) + 2 \sum_{i=2}^n \frac12 + 2 \sum_{i=n+1}^{2n} (e'_i + \frac12 - \beta)
$$
miles in $[\beta,\infty)$. These $\alpha$ miles add up to a walk from position $\beta$ to position $\max\{e'_1,\beta\}$, of which 
$$
\frac{\alpha - (\max\{e'_1,\beta\} - \beta)}{2}
$$
miles are in backward direction and
\begin{eqnarray*}
\alpha' &=& \frac{\alpha + (\max\{e'_1,\beta\} - \beta)}{2} \\
&\le& \frac12\big(e'_1 + \frac12 - \beta\big) + \sum_{i=2}^n \frac12 + \sum_{i=n+1}^{2n} (e'_i + \frac12 - \beta) \\
&=& \frac12 (e'_1 - e'_{2n}) + \sum_{i=n+1}^{2n} e'_i + \frac{2n-1}2 - n \cdot \beta
\end{eqnarray*}
miles are in forward direction.

The camel carries bananas to positions $e'_1, e'_2, \ldots, e'_n$, so 
$$
\alpha' \ge \sum_{i=1}^n (\max\{e'_i,\beta\} - \beta) \ge  \sum_{i=1}^n e'_i - n \cdot \beta
$$
Combining the upper bound and the lower bound of $\alpha'$, we obtain
$$
\sum_{i=1}^n e'_i \le \frac12 (e'_1 - e'_{2n}) + \sum_{i=n+1}^{2n} e'_i + \frac{2n-1}2
$$
Consequently
\begin{eqnarray*}
\frac12 s'_n &\le& - \frac12 e'_{2n} + \sum_{i=n+1}^{2n} e_i + \frac{2n - 1}2 \\
&=& (\frac14 s'_{2n-1} - \frac14 s'_{2n}) + (\frac12 s'_{2n} - \frac12 s'_n) + \frac{2n - 1}2
\end{eqnarray*}
This proves lemma B, and, as a consequence, theorem III. 
\hfill $\Box$

\bigskip \noindent
In some cases it is possible to exchange certain steps in the UWC-strategy
without affecting the final result. Take, e.g., $N = 7$, where the
first two steps may be interchanged. In general, it may be shown that
any optimal strategy $S$ is equivalent to a permutation of the
steps of the UWC-strategy. The proof is left as an exercise to the
reader.

\end{document}